\documentclass[times]{kristall}	        

\usepackage{epsf}
\usepackage{amsmath,amssymb}
\DeclareMathOperator{\norm}{norm}

\newcommand{\CC}{\mathbb{C}}
\newcommand{\ZZ}{\mathbb{Z}}
\newcommand{\QQ}{\mathbb{Q}\hspace{0.5pt}}
\newcommand{\ff}[2]{\genfrac{}{}{0.3pt}{}
{\rule[-0.25ex]{0.3pt}{0.5ex}\makebox[2ex]{$\scriptstyle #2$}}
{\makebox[2ex]{$\scriptstyle #1$}\rule[0.8ex]{0.3pt}{0.5ex}}\hspace{2.5pt}}
\newcommand{\+}{\hfill + \hfill}
\newcommand{\pf}{{\sc Proof}.\quad}
\newcommand{\qed}{\hspace*{\fill}$\square$\medskip}
\newcommand{\I}{\mathrm{i}}
 
\newtheorem{defn}{Definition}
\newtheorem{lem}{Lemma}
\newtheorem{prop}{Proposition}
\newtheorem{thm}{Theorem}

\setbox\Copyrightbox\hbox{\vphantom{Yy}}
\makeatletter
\renewcommand{\@evenhead}{\underline{\hbox to 
\textwidth{{\PageFont\thepage}\hfill\MarkFont\strut \EvenRunningHead}}}%
\renewcommand{\@oddhead}{\underline{\hbox to \textwidth{{\MarkFont 
\strut\OddRunningHead}\hfill\PageFont\thepage}}}%
\renewcommand{\ps@myheadings}{}

\renewcommand{\Pedigree}{%
\begin{tabular}[t]{@{} l @{\hspace{1mm}} l @{}}
\multicolumn{2}{@{} l @{}}{\AuthorsFont\AuthorsComplete\vrule height0pt 
depth 4.83mm width0pt}\\
\raisebox{0pt}[0pt][0pt]{$^{\mathrm I}$} & \InstitutionOne\\
$^{\mathrm {II}}$ &\InstitutionTwo\vrule height0pt depth 4.83mm width0pt\strut
\end{tabular}}

\makeatother

\begin{document}
       
\maketitle		

\begin{abstract}	
\Abstract
\end{abstract}		

\bigskip
\begin{small}
\noindent
{\bf Key Words}:\hspace{\parindent}
Colourings, Planar Modules, Cyclotomic Fields,
Dirichlet Series
\end{small}

\section{Introduction}

Colour symmetries of crystals \cite{Marj1,Marj2,S1,S2,R} and, more
recently, of quasicrystals \cite{MP,L1,B} continue to attract a lot of
attention, simply because so little is known about their
classification, see \cite{L1} for a recent review.  A first step in
this analysis consists in answering the question of how many different
colourings of an infinite point set exist which are compatible with
its underlying symmetry. More precisely, one has to determine the
possible numbers of colours, and to count the corresponding
possibilities to distribute the colours over the point set (up to
permutations), in line with all compatible symmetry constraints.

In this generality, the problem has not even been solved for simple
lattices. One common restriction is to demand that one colour occupies
a subset which is of the same Bravais type as the original set, while
the other colours encode the cosets. Of particular interest are the
cases where the point symmetry is irreducible.  In this situation, to
which we shall also restrict ourselves, several results are known and
can be given in closed form \cite{B,BHLS,BM,L1,L2,R}.

Particularly interesting are planar cases because, on the one hand,
they show up in quasicrystalline $T$-phases, and, on the other hand,
they are linked to the rather interesting classification of planar
Bravais classes with $n$-fold symmetry \cite{MRW}, which is based on a
connection to algebraic number theory in general, and to cyclotomic
fields in particular. The Bravais types correspond to ideal classes
and are unique for the following 29 choices of $n$,
\begin{equation}
\begin{array}{lcl}
n&\in&\{3,4,5,7,8,9,11,12,13,15,16,\\
  &&\hphantom{\{}17,19,20,21,24,25,27,28,32,\\
  &&\hphantom{\{}33,35,36,40,44,45,48,60,84\}.
\end{array}
\label{liste}
\end{equation}
The canonical representatives are the sets of cyclotomic integers
${\mathcal{M}}_{n}=\mathbb{Z}[\xi_{n}]$, the ring of polynomials in
$\xi_{n}$, a primitive $n$th root of unity. To be explicit (which is
not necessary), we choose $\xi_{n}=\exp(2\pi\I/n)$. Apart from $n=1$
(where ${\mathcal{M}}_{1}=\mathbb{Z}$ is one-dimensional), the values
of $n$ in \eqref{liste} correspond to all cases where
$\mathbb{Z}[\xi_{n}]$ is a principal ideal domain and thus has class
number one, see \cite{W,BM} for details. If $n$ is odd, we have
${\mathcal{M}}_{n} = {\mathcal{M}}_{2n}$.  Consequently,
${\mathcal{M}}_{n}$ has $N$-fold symmetry where
\begin{equation}
N\; =\; N(n)\;=\;
    \begin{cases}
         n  & \text{if $n$ is even,} \\
         2n & \text{if $n$ is odd.}
    \end{cases}
\label{symm}
\end{equation}
To avoid duplication of results, values $n\equiv 2\bmod 4$ do not
appear in the above list \eqref{liste}. There are systematic
mathematical reasons to prefer this convention to the notation of
\cite{L1,MRW}, some of which will appear later on.

It is this very connection to algebraic number theory which results in
the Bravais classification \cite{MRW}, and also allows for a solution
of the combinatorial part of the colouring problem by means of
Dirichlet series generating functions, compare \cite{B,BM}. The values
of $n$ from our list \eqref{liste} are naturally grouped according to
$\phi(n)$, which is Euler's totient function
\begin{equation} \label{eulerphi}
   \phi(n) \; = \; \bigl\lvert\{ 1\le k \le n \mid \gcd(k,n)=1 \}
    \bigr\rvert\,.
\end{equation}
Note that $\phi(n)=2$ covers the two crystallographic cases $n=3$
(triangular lattice) and $n=4$ (square lattice), while
\mbox{$\phi(n)=4$} means $n\in\{5,8,12\}$ which are the standard
symmetries of genuine planar quasicrystals. Again, $n=10$ is covered
implicitly, as explained above.

The methods emerging from the connection to number theory are useful
for the description of aperiodic order in general \cite{P}. In the
planar case, the link to cyclotomic fields allows the full treatment
of all 29 cases listed above, which was observed and used in several
papers \cite{MRW,PBR,B,BM}.  Here, we shall explain the use of
cyclotomic fields more explicitly, and present detailed results for
the combinatorial part of the colouring problem in all 29 cases. We
shall not follow the standard approach of algebraic number theory (see
\cite{BS} for a good and readable general introduction), but take
the slightly unusual point of view of Dirichlet series generating
functions. This allows for a rather straightforward calculation of all
quantities required, starting from the representation theory of finite
Abelian groups. Furthermore, we shall also present some asymptotic
results that can be extracted from the generating function.

The modules arise in quasicrystal theory in several ways. One is
through the Fourier module that supports the Bragg peaks. Another,
more important, is via the limit translation module of a (discrete)
quasiperiodic tiling. Characteristic points of the latter (e.g.,
vertex points) are then model sets (or cut and project sets) on the
basis of the entire (dense) module, seen as a lattice in a space of
higher dimension, see \cite{B,BGS1,BHLS} for details. This way, one
obtains a one-to-one relation between colourings of discrete objects
and of the underlying dense module. Since the latter is universal, it
is the natural object to study.

Previously, results for the cases with $\phi(n)\le 10$ were given in
\cite{B,L1,BGS1,BGS2}, mainly without proofs. Other combinatorial
problems of (quasi)crystallography have been addressed by similar
methods, compare \cite{BG}.

\section{Number theoretic formulation}

In view of the above remarks, we formulate the problem immediately in
terms of the full modules.  Consider ${\mathcal{M}}_n$ with a fixed
$n$, e.g., from our list~\eqref{liste}.  Then, ${\mathcal{M}}_n$ is an
Abelian group, and also a $\mathbb{Z}$-module of rank $\phi(n)$. We
view it as a subset of the Euclidean plane, identified with $\CC$, and
hence as a geometric object. In this setting, important subgroups of
${\mathcal{M}}_n$ are those that are images of ${\mathcal{M}}_n$ under
a similarity transformation, i.e., a mapping
\[
z\mapsto az \quad\mbox{or}\quad
z\mapsto a\bar{z}\qquad
\mbox{with~~$a\in\CC^{*}:=\CC\setminus\{0\}$}\,.
\]
Such subgroups are called {\em similarity submodules}, see \cite{BM}
for details.  Similarity submodules are very natural objects
algebraically, and have been studied in the context of colouring
problems in dimensions two, three and four \cite{B,BM,BM2}. To make
this more precise, let us start with a definition.
 
\begin{defn}\label{def:bravais}
A {\em Bravais colouring} of the module ${\mathcal{M}}_{n}$ with $k$
colours is a mapping\/ $c\!: {\mathcal{M}}_{n}\rightarrow
\{1,2,\ldots,k\}$ such that one of the sets $c^{-1}(\ell)$, for
$1\le\ell\le k$, is a similarity submodule of ${\mathcal{M}}_{n}$ of
index $k$, and the others, one by one, are the corresponding cosets in
${\mathcal{M}}_{n}$. The number of Bravais colourings of
${\mathcal{M}}_{n}$ with $k$ colours, up to permutations of the
colours, is denoted by\/ $a^{}_{n}(k)$.
\end{defn}

In other words, a Bravais colouring $c$ of ${\mathcal{M}}_{n}$ is a
partition
\begin{equation}\label{part}
{\mathcal{M}}_{n}\;=\;\dot{\bigcup_{1\le\ell\le k}} c^{-1}(\ell)
\end{equation}
into $k$ disjoint sets that are translates of each other and carry
pairwise different colours. In addition, one of the sets is a
similarity submodule of ${\mathcal{M}}_{n}$. The interest in this kind
of partition originates from its rather important group-theoretical
structure, compare \cite{S1,S2,Marj1,Marj2,L1} and references
therein. Also, each such partition induces a colouring of tilings
that are compatible with the module (usually formulated in terms of local
indistinguishability, compare \cite{L1}), see \cite{BGS1,BHLS,SL} for
examples.

Here, we concentrate on the combinatorial aspect, i.e., on
the determination of the values of $a^{}_{n}(k)$, which define an
integer-valued arithmetic function.

\begin{table*}[t]
\caption[]{Basic indices for the unramified primes of
$\ZZ[\xi_{n}]$ with $n$ from \eqref{liste}. 
The symbol $\protect\ff{k}{\ell}$ means that primes 
$p\equiv k\bmod n$ contribute via $p^{\ell}$ as basic 
index, where $\ell$ is the smallest integer such that
$k^{\ell}\equiv 1\bmod n$, and the integer $m$ entering Eq.~\eqref{factor} is 
$m=\phi(n)/\ell$.\label{tab:basic}}
\def\Strut{\large\strut}
\renewcommand{\arraystretch}{1.25}
\begin{tabular*}{\textwidth}{@{\extracolsep{\fill}} @{}crl@{}}  
 \hline
\vspace{-0.5pt}$\phi(n)$ & $n$ & general primes $p$ \Strut\\
 \hline
 2 &  3 & $\ff{1}{1}$ $\ff{2}{2}$ \Strut\\
   &  4 & $\ff{1}{1}$ $\ff{3}{2}$ \rule[-1.5ex]{0ex}{1.5ex}\Strut\\
 \hline
 4 &  5 & $\ff{1}{1}$ $\ff{2}{4}$ 
          $\ff{3}{4}$ $\ff{4}{2}$\Strut\\
   &  8 & $\ff{1}{1}$ $\ff{3}{2}$ 
          $\ff{5}{2}$ $\ff{7}{2}$\Strut\\
   & 12 & $\ff{1}{1}$ $\ff{5}{2}$ 
          $\ff{7}{2}$ $\ff{11}{2}$\rule[-1.5ex]{0ex}{1.5ex}\Strut\\
 \hline
 6 &  7 & $\ff{1}{1}$ $\ff{2}{3}$ 
          $\ff{3}{6}$ $\ff{4}{3}$ 
          $\ff{5}{6}$ $\ff{6}{2}$\Strut\\
   &  9 & $\ff{1}{1}$ $\ff{2}{6}$ 
          $\ff{4}{3}$ $\ff{5}{6}$ 
          $\ff{7}{3}$ $\ff{8}{2}$\rule[-1.5ex]{0ex}{1.5ex}\Strut\\
 \hline
 8 & 15 & $\ff{1}{1}$ $\ff{2}{4}$ $\ff{4}{2}$ $\ff{7}{4}$ 
          $\ff{8}{4}$ $\ff{11}{2}$ $\ff{13}{4}$ $\ff{14}{2}$\Strut\\
   & 16 & $\ff{1}{1}$ $\ff{3}{4}$ $\ff{5}{4}$ $\ff{7}{2}$ 
          $\ff{9}{2}$ $\ff{11}{4}$ $\ff{13}{4}$ $\ff{15}{2}$\Strut\\
   & 20 & $\ff{1}{1}$ $\ff{3}{4}$ $\ff{7}{4}$ $\ff{9}{2}$ 
          $\ff{11}{2}$ $\ff{13}{4}$ $\ff{17}{4}$ $\ff{19}{2}$\Strut\\
   & 24 & $\ff{1}{1}$ $\ff{5}{2}$ $\ff{7}{2}$ $\ff{11}{2}$ 
          $\ff{13}{2}$ $\ff{17}{2}$ $\ff{19}{2}$ $\ff{23}{2}$
\rule[-1.5ex]{0ex}{1.5ex}\Strut\\
\hline
10 & 11 & $\ff{1}{1}$ $\ff{2}{10}$ $\ff{3}{5}$ $\ff{4}{5}$
          $\ff{5}{5}$ $\ff{6}{10}$ $\ff{7}{10}$ $\ff{8}{10}$
          $\ff{9}{5}$ $\ff{10}{2}$\rule[-1.5ex]{0ex}{1.5ex}\Strut\\ 
\hline
12 & 13 & $\ff{1}{1}$ $\ff{2}{12}$ $\ff{3}{3}$ $\ff{4}{6}$
          $\ff{5}{4}$ $\ff{6}{12}$ $\ff{7}{12}$ $\ff{8}{4}$
          $\ff{9}{3}$ $\ff{10}{6}$ $\ff{11}{12}$ $\ff{12}{2}$\Strut\\
   & 21 & $\ff{1}{1}$ $\ff{2}{6}$ $\ff{4}{3}$ $\ff{5}{6}$
          $\ff{8}{2}$ $\ff{10}{6}$ $\ff{11}{6}$ $\ff{13}{2}$
          $\ff{16}{3}$ $\ff{17}{6}$ $\ff{19}{6}$ $\ff{20}{2}$\Strut\\
   & 28 & $\ff{1}{1}$ $\ff{3}{6}$ $\ff{5}{6}$ $\ff{9}{3}$
          $\ff{11}{6}$ $\ff{13}{2}$ $\ff{15}{2}$ $\ff{17}{6}$
          $\ff{19}{6}$ $\ff{23}{6}$ $\ff{25}{3}$ $\ff{27}{2}$\Strut\\
   & 36 & $\ff{1}{1}$ $\ff{5}{6}$ $\ff{7}{6}$ $\ff{11}{6}$
          $\ff{13}{3}$ $\ff{17}{2}$ $\ff{19}{2}$ $\ff{23}{6}$
          $\ff{25}{3}$ $\ff{29}{6}$ $\ff{31}{6}$ $\ff{35}{2}$
\rule[-1.5ex]{0ex}{1.5ex}\Strut\\
\hline
16 & 17 & $\ff{1}{1}$ $\ff{2}{8}$ $\ff{3}{16}$ $\ff{4}{4}$
          $\ff{5}{16}$ $\ff{6}{16}$ $\ff{7}{16}$ $\ff{8}{8}$
          $\ff{9}{8}$ $\ff{10}{16}$ $\ff{11}{16}$ $\ff{12}{16}$
          $\ff{13}{4}$ $\ff{14}{16}$ $\ff{15}{8}$ $\ff{16}{2}$\Strut\\
   & 32 & $\ff{1}{1}$ $\ff{3}{8}$ $\ff{5}{8}$ $\ff{7}{4}$
          $\ff{9}{4}$ $\ff{11}{8}$ $\ff{13}{8}$ $\ff{15}{2}$
          $\ff{17}{2}$ $\ff{19}{8}$ $\ff{21}{8}$ $\ff{23}{4}$
          $\ff{25}{4}$ $\ff{27}{8}$ $\ff{29}{8}$ $\ff{31}{2}$\Strut\\
   & 40 & $\ff{1}{1}$ $\ff{3}{4}$ $\ff{7}{4}$ $\ff{9}{2}$
          $\ff{11}{2}$ $\ff{13}{4}$ $\ff{17}{4}$ $\ff{19}{2}$
          $\ff{21}{2}$ $\ff{23}{4}$ $\ff{27}{4}$ $\ff{29}{2}$
          $\ff{31}{2}$ $\ff{33}{4}$ $\ff{37}{4}$ $\ff{39}{2}$\Strut\\
   & 48 & $\ff{1}{1}$ $\ff{5}{4}$ $\ff{7}{2}$ $\ff{11}{4}$
          $\ff{13}{4}$ $\ff{17}{2}$ $\ff{19}{4}$ $\ff{23}{2}$
          $\ff{25}{2}$ $\ff{29}{4}$ $\ff{31}{2}$ $\ff{35}{4}$
          $\ff{37}{4}$ $\ff{41}{2}$ $\ff{43}{4}$ $\ff{47}{2}$\Strut\\
   & 60 & $\ff{1}{1}$ $\ff{7}{4}$ $\ff{11}{2}$ $\ff{13}{4}$
          $\ff{17}{4}$ $\ff{19}{2}$ $\ff{23}{4}$ $\ff{29}{2}$
          $\ff{31}{2}$ $\ff{37}{4}$ $\ff{41}{2}$ $\ff{43}{4}$
          $\ff{47}{4}$ $\ff{49}{2}$ $\ff{53}{4}$ $\ff{59}{2}$
\rule[-1.5ex]{0ex}{1.5ex}\Strut\\
\hline
18 & 19 & $\ff{1}{1}$ $\ff{2}{18}$ $\ff{3}{18}$ $\ff{4}{9}$
          $\ff{5}{9}$ $\ff{6}{9}$ $\ff{7}{3}$ $\ff{8}{6}$
          $\ff{9}{9}$ $\ff{10}{18}$ $\ff{11}{3}$ $\ff{12}{6}$
          $\ff{13}{18}$ $\ff{14}{18}$ $\ff{15}{18}$ $\ff{16}{9}$
          $\ff{17}{9}$ $\ff{18}{2}$\Strut\\
   & 27 & $\ff{1}{1}$ $\ff{2}{18}$ $\ff{4}{9}$ $\ff{5}{18}$
          $\ff{7}{9}$ $\ff{8}{6}$ $\ff{10}{3}$ $\ff{11}{18}$
          $\ff{13}{9}$ $\ff{14}{18}$ $\ff{16}{9}$ $\ff{17}{6}$
          $\ff{19}{3}$ $\ff{20}{18}$ $\ff{22}{9}$ $\ff{23}{18}$
          $\ff{25}{9}$ $\ff{26}{2}$
\rule[-1.5ex]{0ex}{1.5ex}\Strut\\ 
\hline
20 & 25 & $\ff{1}{1}$ $\ff{2}{20}$ $\ff{3}{20}$ $\ff{4}{10}$
          $\ff{6}{5}$ $\ff{7}{4}$ $\ff{8}{20}$ $\ff{9}{10}$
          $\ff{11}{5}$ $\ff{12}{20}$ $\ff{13}{20}$ $\ff{14}{10}$
          $\ff{16}{5}$ $\ff{17}{20}$ $\ff{18}{4}$ $\ff{19}{10}$
          $\ff{21}{5}$ $\ff{22}{20}$ $\ff{23}{20}$ $\ff{24}{2}$\Strut\\
   & 33 & $\ff{1}{1}$ $\ff{2}{10}$ $\ff{4}{5}$ $\ff{5}{10}$
          $\ff{7}{10}$ $\ff{8}{10}$ $\ff{10}{2}$ $\ff{13}{10}$
          $\ff{14}{10}$ $\ff{16}{5}$ $\ff{17}{10}$ $\ff{19}{10}$
          $\ff{20}{10}$ $\ff{23}{2}$ $\ff{25}{5}$ $\ff{26}{10}$
          $\ff{28}{10}$ $\ff{29}{10}$ $\ff{31}{5}$ $\ff{32}{2}$\Strut\\ 
   & 44 & $\ff{1}{1}$ $\ff{3}{10}$ $\ff{5}{5}$ $\ff{7}{10}$
          $\ff{9}{5}$ $\ff{13}{10}$ $\ff{15}{10}$ $\ff{17}{10}$
          $\ff{19}{10}$ $\ff{21}{2}$ $\ff{23}{2}$ $\ff{25}{5}$
          $\ff{27}{10}$ $\ff{29}{10}$ $\ff{31}{10}$ $\ff{35}{10}$
          $\ff{37}{5}$ $\ff{39}{10}$ $\ff{41}{10}$ $\ff{43}{2}$
\rule[-1.5ex]{0ex}{1.5ex}\Strut\\
\hline
24 & 35 & $\ff{1}{1}$ $\ff{2}{12}$ $\ff{3}{12}$ $\ff{4}{6}$
          $\ff{6}{2}$ $\ff{8}{4}$ $\ff{9}{6}$ $\ff{11}{3}$
          $\ff{12}{12}$ $\ff{13}{4}$ $\ff{16}{3}$ $\ff{17}{12}$
          $\ff{18}{12}$ $\ff{19}{6}$ $\ff{22}{4}$ $\ff{23}{12}$
          $\ff{24}{6}$ $\ff{26}{6}$ $\ff{27}{4}$ $\ff{29}{2}$
          $\ff{31}{6}$ $\ff{32}{12}$ $\ff{33}{12}$ $\ff{34}{2}$\Strut\\
   & 45 & $\ff{1}{1}$ $\ff{2}{12}$ $\ff{4}{6}$ $\ff{7}{12}$
          $\ff{8}{4}$ $\ff{11}{6}$ $\ff{13}{12}$ $\ff{14}{6}$
          $\ff{16}{3}$ $\ff{17}{4}$ $\ff{19}{2}$ $\ff{22}{12}$
          $\ff{23}{12}$ $\ff{26}{2}$ $\ff{28}{4}$ $\ff{29}{6}$
          $\ff{31}{3}$ $\ff{32}{12}$ $\ff{34}{6}$ $\ff{37}{4}$
          $\ff{38}{12}$ $\ff{41}{6}$ $\ff{43}{12}$ $\ff{44}{2}$\Strut\\
   & 84 & $\ff{1}{1}$ $\ff{5}{6}$ $\ff{11}{6}$ $\ff{13}{2}$
          $\ff{17}{6}$ $\ff{19}{6}$ $\ff{23}{6}$ $\ff{25}{3}$
          $\ff{29}{2}$ $\ff{31}{6}$ $\ff{37}{3}$ $\ff{41}{2}$
          $\ff{43}{2}$ $\ff{47}{6}$ $\ff{53}{6}$ $\ff{55}{2}$
          $\ff{59}{6}$ $\ff{61}{6}$ $\ff{65}{6}$ $\ff{67}{6}$
          $\ff{71}{2}$ $\ff{73}{6}$ $\ff{79}{6}$ $\ff{83}{2}$
\rule[-1.5ex]{0ex}{1.5ex}\Strut\\
\hline
\end{tabular*}
\renewcommand{\arraystretch}{1}
\end{table*}

\begin{lem}\label{lem:ideal}
Let\/ ${\mathcal{M}}_{n}=\ZZ[\xi_{n}]$ with\/ $n\ge 3$. A similarity
submodule of\/ ${\mathcal{M}}_{n}$ is a principal ideal of\/
$\ZZ[\xi_{n}]$. Consequently, $a^{}_{n}(k)$ is the number of principal
ideals of\/ $\ZZ[\xi_{n}]$ of norm\/ $k$.
\end{lem}
\pf
Let us observe that the module $\ZZ[\xi_{n}]$ is invariant under
complex conjugation.  Consequently, all similarity submodules of
$\ZZ[\xi_{n}]$ are of the form $a\ZZ[\xi_{n}]$ with $a\in\CC^{*}$. As
\mbox{$a\ZZ[\xi_{n}]\subset\ZZ[\xi_{n}]$}, we must have $0\ne
a\in\ZZ[\xi_{n}]$.  So, a similarity submodule is a principal ideal,
and vice versa. Each principal ideal gives rise to precisely one
Bravais colouring $c$ of $\ZZ[\xi_{n}]$, up to permutations of the
colours, where, without loss of generality, $c^{-1}(1)$ is the
principal ideal. The norm of an ideal equals the subgroup index
$k=\big[\ZZ[\xi_{n}]:a\ZZ[\xi_{n}]\big]$, hence counts the number of
cosets, resp.\ colours. In view of Eq.~\eqref{part}, and the meaning
of the sets $c^{-1}(\ell)$, our assertion follows.
\qed

For $n$ from the list~\eqref{liste}, {\em all}\/ ideals of
$\ZZ[\xi_{n}]$ are principal, because the corresponding cyclotomic
fields $\QQ(\xi_{n})$ have class number one
\cite[Thm.~11.1]{W}. This list is exhaustive except for the case $n=1$
which corresponds to the field $\QQ$ itself. Our
combinatorial problem then amounts to counting all (non-zero) ideals
of a given index. Consequently, the arithmetic function $a^{}_{n}(k)$
is {\em multiplicative}, i.e.,
$a^{}_{n}(k\ell)=a^{}_{n}(k)a^{}_{n}(\ell)$ for $k,\ell$ coprime, and
$a^{}_{n}(1)=1$. A succinct way to encapsulate the numbers
$a^{}_{n}(k)$ is the use of a Dirichlet series generating function
\cite{Wilf},
\begin{equation}\label{def:diri}
 F_{n} (s) \; := \;
 \sum_{k=1}^{\infty} \frac{a^{}_{n}(k)}{k^s}\, .
\end{equation}
As we shall see, this function converges for all $s\in\CC$ with ${\rm
Re}(s)>1$.  The connection to the cyclotomic fields permits the direct
calculation of this Dirichlet series by means of algebraic number
theory.
\begin{prop}\label{prop:dede}
For all\/ $n$ in the list\/ $\eqref{liste}$, the Dirichlet series\/
$F_{n}(s)$ of Eq.~$\eqref{def:diri}$ equals the Dedekind zeta function
of the cyclotomic field\/ $\QQ(\xi_{n})$, i.e.,
\[
F_{n}(s)\;=\;\zeta^{}_{\QQ(\xi_{n})}(s)\;:=\;
\sum_{\mathfrak{a}}\,\frac{1}{\norm(\mathfrak{a})^{s}}\, ,
\quad {\rm Re}(s)>1\, ,
\]
where\/ $\mathfrak{a}$ runs over all non-zero ideals of\/
$\ZZ[\xi_{n}]$.
\end{prop}
\pf
The Dedekind zeta function of $\QQ(\xi_{n})$ is, by definition, the
Dirichlet series generating function of the number of ideals of a
given index in the maximal order $\ZZ[\xi_{n}]$. Since $n$ is from our
list~\eqref{liste}, we are in the class number one case, hence all
ideals are principal, and the claim follows from
Lemma~\ref{lem:ideal}. The convergence result is standard
\cite{W}.
\qed

For general $n$, the arithmetic function $a_{n}^{}(k)$ still counts
the principal ideals of $\QQ(\xi_{n})$. Whenever non-principal ideals
exist, $F_{n}(s)$ is no longer given by the Dedekind zeta function of
$\QQ(\xi_{n})$. For the related problem of coincidence site modules,
an example ($n=23$) is treated in
\cite{PBR}.

\begin{table}
\caption[]{List of all ramified primes with corresponding
integers $\ell$ and $m$ for $\ZZ[\xi_{n}]$ with $n$ from the list
\eqref{liste}. Here, $r$ is the $p$-free part of $n$, and $\ell\, m
=\phi(r)$. Details on the connection to Table~\ref{tab:basic} are
given in the text. \label{tab:rami}}
\def\Strut{\large\strut}
\small
\begin{tabular*}{\linewidth}{@{\extracolsep{\fill}} @{} rrrrrrr @{}}
\hline
$\phi(n)$ &
$n$ &
$p$ &
$r$ &
$\phi(r)$ &
$\ell$ &
$m$ \Strut\\
\hline
 2 &  3 &  3 & 1 & 1 &  1 & 1 \Strut\\ 
   &  4 &  2 & 1 & 1 &  1 & 1 \Strut\\ \hline
 4 &  5 &  5 & 1 & 1 &  1 & 1 \Strut\\ 
   &  8 &  2 & 1 & 1 &  1 & 1 \Strut\\ 
   & 12 &  2 & 3 & 2 &  2 & 1 \Strut\\ 
   &    &  3 & 4 & 2 &  2 & 1 \Strut\\ \hline
 6 &  7 &  7 & 1 & 1 &  1 & 1 \Strut\\ 
   &  9 &  3 & 1 & 1 &  1 & 1 \Strut\\ \hline
 8 & 15 &  3 & 5 & 4 &  4 & 1 \Strut\\ 
   &    &  5 & 3 & 2 &  2 & 1 \Strut\\ 
   & 16 &  2 & 1 & 1 &  1 & 1 \Strut\\ 
   & 20 &  2 & 5 & 4 &  4 & 1 \Strut\\ 
   &    &  5 & 4 & 2 &  1 & 2 \Strut\\ 
   & 24 &  2 & 3 & 2 &  2 & 1 \Strut\\ 
   &    &  3 & 8 & 4 &  2 & 2 \Strut\\ \hline
10 & 11 & 11 & 1 & 1 &  1 & 1 \Strut\\ \hline
12 & 13 & 13 & 1 & 1 &  1 & 1 \Strut\\ 
   & 21 &  3 & 7 & 6 &  6 & 1 \Strut\\ 
   &    &  7 & 3 & 2 &  1 & 2 \Strut\\ 
   & 28 &  2 & 7 & 6 &  3 & 2 \Strut\\ 
   &    &  7 & 4 & 2 &  2 & 1 \Strut\\ 
   & 36 &  2 & 9 & 6 &  6 & 1 \Strut\\ 
   &    &  3 & 4 & 2 &  2 & 1 \Strut\\ \hline
16 & 17 & 17 & 1 & 1 &  1 & 1\Strut\\ 
   & 32 &  2 & 1 & 1 &  1 & 1\Strut\\ 
   & 40 &  2 & 5 & 4 &  4 & 1\Strut\\
   &    &  5 & 8 & 4 &  2 & 2\Strut\\
   & 48 &  2 & 3 & 2 &  2 & 1\Strut\\
   &    &  3 &16 & 8 &  4 & 2\Strut\\ 
   & 60 &  2 &15 & 8 &  4 & 2\Strut\\
   &    &  3 &20 & 8 &  4 & 2\Strut\\ 
   &    &  5 &12 & 4 &  2 & 2\Strut\\ \hline
18 & 19 & 19 & 1 & 1 &  1 & 1\Strut\\ 
   & 27 &  3 & 1 & 1 &  1 & 1\Strut\\ \hline
20 & 25 &  5 & 1 & 1 &  1 & 1\Strut\\
   & 33 &  3 &11 &10 &  5 & 2\Strut\\
   &    & 11 & 3 & 2 &  2 & 1\Strut\\
   & 44 &  2 &11 &10 & 10 & 1\Strut\\ 
   &    & 11 & 4 & 2 &  2 & 1\Strut\\ \hline
24 & 35 &  5 & 7 & 6 &  6 & 1\Strut\\
   &    &  7 & 5 & 4 &  4 & 1\Strut\\
   & 45 &  3 & 5 & 4 &  4 & 1\Strut\\
   &    &  5 & 9 & 6 &  6 & 1\Strut\\
   & 84 &  2 &21 &12 &  6 & 2\Strut\\
   &    &  3 &28 &12 &  6 & 2\Strut\\
   &    &  7 &12 & 4 &  2 & 2\Strut\\
\hline
\end{tabular*}
\end{table}

The multiplicativity of the arithmetic function $a^{}_{n}(k)$ for $n$
in the list~\eqref{liste} implies that it is sufficient to know the
values $a^{}_{n}(p^r)$ for all primes $p$ and powers $r>0$. On the
level of the generating function, this corresponds to the Euler
product expansion of $F_{n}(s)$, again for $\{{\rm Re}(s)>1\}$,
\begin{equation}
 F_{n}(s) \; := \;
 \sum_{k=1}^{\infty} \frac{a^{}_{n}(k)}{k^s} \; = \;
 \prod_{p} E_{n}(p^{-s})
 \label{euler}
\end{equation}
where $p$ runs over the primes of $\ZZ$. 

As we shall see below, each Euler factor is of the form
\begin{eqnarray}
 E_{n}(p^{-s}) & = & \frac{1}{(1-p^{-\ell s})^m}\nonumber\\
 & = &
 \sum_{j=0}^{\infty} \binom{j+m-1}{m-1}
 \frac{1}{(p^{s})^{\ell j}}
 \label{factor}
\end{eqnarray}
from which one extracts the values of $a^{}_{n}(p^r)$ for $r\ge 0$.
The integers $\ell$ and $m$ are characteristic quantities that depend
on $n$ and on the prime $p$. In particular, $m$ is the number of prime
ideal divisors of $p$, and $\ell$ is their degree (also called the
residue class degree) \cite{W}. With the usual approach of algebraic
number theory, its determination (and that of the corresponding $m$)
depends on whether $p$ divides $n$ or not. This can be circumvented
with a formula based on Dirichlet characters, as we shall explain
below.

If $p$ and $n$ are coprime, we have $p\equiv k \bmod n$ for some $k$
such that $k$ and $n$ are coprime, i.e., $\gcd(k,n)=1$.  In this case,
$\ell\, m =\phi(n)$, which fixes $m$. Furthermore, the residue class
degree $\ell$ is the smallest integer such that $k^{\ell}\equiv 1\bmod
n$, see \cite[Thm.~2.3]{W}. The cases of these primes are listed as
$\ff{k}{\ell}$ in Table~\ref{tab:basic}.

In addition, for each $n$, there are finitely many primes $p$ which
divide $n$, the so-called {\em ramified}\/ primes. In this case,
$n=r\,p^t$ with $t\ge 1$ and $r$ an integer not divisible by $p$, so
that $r$ is the $p$-free part of $n$. The values of $\ell$ and $m$
needed in the case $p|n$ equal those needed for the situation where
$n$ is replaced by $r$, which brings us back to the previous case (the
proof is more involved and can be extracted from \cite{L}). In
particular, $\ell\, m = \phi(r)$, and $\ell=m=1$ whenever $r=1$.  The
complete result is listed in Table \ref{tab:rami} for
convenience. With this information, one can easily calculate the
possible numbers of colours and the generating functions by inserting
\eqref{factor} into \eqref{euler} and expanding the Euler product,
which is an easy task for an algebraic program package.

Let us briefly come back to Eq.~\eqref{def:diri}. If $a_{n}(k)$ is
multiplicative, the set of possible numbers of colours forms a
semigroup with unit. It is generated by the basic indices which are
obtained as $p^{\ell}$ from Tables \ref{tab:basic} and
\ref{tab:rami}. In more general situations, this semigroup structure
is lost, compare the related coincidence problem \cite{PBR}.

\section{Dirichlet characters and zeta functions}

The purpose of this section is to summarise some basic properties of
cyclotomic fields that are helpful to calculate the Dedekind zeta
functions explicitly, and hence also the characteristic integers
$\ell$ and $m$ via Eq.~\eqref{factor}, in a unified fashion. 
Moreover, some analytic properties of the zeta functions will
become accessible this way, which we shall need to determine
asymptotic properties of our counting functions. The approach
uses the well-known theory of Dirichlet characters, which
we shall now summarise.

The object we start from is the {\em Galois group}\/ $G_n$ of the
(cyclotomic) field extension $\QQ(\xi_{n})/\QQ$. Here, the Galois
group simply consists of all automorphisms $\sigma$ of $\QQ(\xi_{n})$
that fix all rational elements, i.e., $\sigma(x)=x$ for all
$x\in\QQ$. It is an Abelian group of order $\phi(n)$. The following
result is standard \cite[Thm.~2.5]{W}.
 
\begin{prop}\label{prop:galois}
The Galois group\/ $G_n$ of\/ $\QQ(\xi_{n})/\QQ$ is Abelian and of
order\/ $\phi(n)$. It is isomorphic to
\[
\{1\le k\le n\mid \gcd(k,n)=1\}
\]
with multiplication modulo\/ $n$ as the group operation, and\/ $k$
representing the automorphism defined by the map\/
$\xi_{n}\mapsto\xi_{n}^{k}$.  For the cases with class number one, the
group structure and a possible set of generators are given in
Table~$\ref{tab:galois}$.\hspace*{\fill}$\square$
\end{prop}

Next, we need the Dirichlet characters, which can be seen as
extensions of the characters $\chi$ of $G_n$ in the sense of
representation theory of (finite) Abelian groups. So, each $\chi$ is a
group homomorphism from $G_n$ into the unit circle in $\CC$.  With
$G_n$ as given in Proposition~\ref{prop:galois}, $\chi$ is defined on
all elements of $\{1\le k\le n\mid \gcd(k,n)=1\}$ modulo $n$. To
obtain a primitive Dirichlet character \cite[Ch.~3]{W}, of which there
are precisely $\phi(n)$ different ones, we have to extend $\chi$ to
{\em all}\/ integers. Clearly, $\chi$ is periodic with period $n$, so
it suffices to define $\chi$ on the missing integers $\le n$. This has
to be done in such a way that the resulting character, which we still
denote by $\chi$, is totally multiplicative, i.e.,
$\chi(k\ell)=\chi(k)\chi(\ell)$ for {\em all}\/ $k,\ell\ge 1$, and
that the period of the resulting character is minimal. This minimal
period, which is always a divisor of $n$, is called the {\em
conductor}\/ of $\chi$, denoted by $f_{\chi}^{}$.  In this process,
one has $\chi(k)=0$ whenever $\gcd(k,f^{}_{\chi})\ne 1$.  An
illustrative example ($n=20$) is shown in Table~\ref{tab:char}.

\begin{table}
\caption[]{Galois groups $G_n$ and their generators (written
as residue classes $\bmod n$) for the 29 non-trivial cyclotomic 
fields with class number one, see Eq.~\eqref{liste}.\label{tab:galois}}
\def\Strut{\large\strut}
\small
\begin{tabular*}{\linewidth}{@{\extracolsep{\fill}} @{} rrll @{}}
\hline
$n$ &
$\phi(n)$ &
Galois group &
generators\Strut\\
\hline
 3 &  2 & $C_2$                             &  (2) \Strut\\ 
 4 &  2 & $C_2$                             &  (3) \Strut\\
 5 &  4 & $C_4$                             &  (2) \Strut\\
 7 &  6 & $C_6$                             &  (3) \Strut\\
 8 &  4 & $C_2\!\times\! C_2$               &  (3),(5) \Strut\\
 9 &  6 & $C_6$                             &  (2) \Strut\\
11 & 10 & $C_{10}$                          &  (2) \Strut\\
12 &  4 & $C_2\!\times\! C_2$               &  (5),(7) \Strut\\
13 & 12 & $C_{12}$                          &  (2) \Strut\\
15 &  8 & $C_4\!\times\! C_2$               &  (2),(11) \Strut\\
16 &  8 & $C_4\!\times\! C_2$               &  (3),(7) \Strut\\
17 & 16 & $C_{16}$                          &  (3) \Strut\\
19 & 18 & $C_{18}$                          &  (2) \Strut\\
20 &  8 & $C_4\!\times\! C_2$               &  (3),(11) \Strut\\
21 & 12 & $C_6\!\times\! C_2$               &  (2),(13) \Strut\\
24 &  8 & $C_2\!\times\! C_2\!\times\! C_2$ &  (5),(7),(13) \Strut\\
25 & 20 & $C_{20}$                          &  (2) \Strut\\  
27 & 18 & $C_{18}$                          &  (2) \Strut\\ 
28 & 12 & $C_6\!\times\! C_2$               &  (3),(13) \Strut\\
32 & 16 & $C_8\!\times\! C_2$               &  (3),(15) \Strut\\
33 & 20 & $C_{10}\!\times\! C_2$            &  (2),(10) \Strut\\
35 & 24 & $C_{12}\!\times\! C_2$            &  (2),(6) \Strut\\
36 & 12 & $C_6\!\times\! C_2$               &  (5),(19) \Strut\\
40 & 16 & $C_4\!\times\! C_2\!\times\! C_2$ &  (3),(11),(21) \Strut\\
44 & 20 & $C_{10}\!\times\! C_2$            &  (3),(21) \Strut\\
45 & 24 & $C_{12}\!\times\! C_2$            &  (2),(26) \Strut\\
48 & 16 & $C_4\!\times\! C_2\!\times\! C_2$ &  (5),(17),(23) \Strut\\
60 & 16 & $C_4\!\times\! C_2\!\times\! C_2$ &  (7),(11),(19) \Strut\\
84 & 24 & $C_6\!\times\! C_2\!\times\! C_2$ &  (5),(13),(43) \Strut\\
\hline
\end{tabular*}
\end{table}

\begin{table*}
\caption[]{Construction of all primitive Dirichlet characters for $n=20$.
The generators of $G_{20}\simeq C_{4}\times C_{2}$ are $g=(3)$ and
$h=(11)$. Each character $\chi$ originates from a product of
characters of $C_{4}$ and $C_{2}$, and is thus labelled by a pair
$(i,j)$. The two extra lines show the characteristic indices $\ell$ and
$m$ for this example.\label{tab:char}}
\def\Strut{\large\strut}
\small
\begin{tabular*}{\textwidth}{@{\extracolsep{\fill}} @{} 
cr|cccccccccccccccccccc @{}}
\hline
       &               & 1 &  2 &  3 &  4 &  5 &  6 &  7 &  8 &  9 & 10 & 
                        11 & 12 & 13 & 14 & 15 & 16 & 17 & 18 & 19 & 20 
\Strut\\
$\chi$ & $f_{\chi}^{}$ &$e$&    & $g$&    &    &    &$g^3$&   &$g^2$&   &
                        $h$&    &$gh$&    &    &   &$g^3h$&  &$g^2h$&   
\Strut\\
\hline
$(0,0)$& 1             & 1 &  1 &  1 &  1 &  1 &  1 &  1 &  1 &  1 &  1 &
                         1 &  1 &  1 &  1 &  1 &  1 &  1 &  1 &  1 &  1 
\Strut\\
$(1,0)$& 5             & 1 &$-\I$& $\I$&$-1$&  0 &  1 &$-\I$& $\I$&$-1$&  0 &
                         1 &$-\I$& $\I$&$-1$&  0 &  1 &$-\I$& $\I$&$-1$&  0 
\Strut\\
$(2,0)$& 5             & 1 &$-1$&$-1$&  1 &  0 &  1 &$-1$&$-1$&  1 &  0 &
                         1 &$-1$&$-1$&  1 &  0 &  1 &$-1$&$-1$&  1 &  0 
\Strut\\
$(3,0)$& 5             & 1 & $\I$&$-\I$&$-1$&  0 &  1 & $\I$&$-\I$&$-1$&  0 &
                         1 & $\I$&$-\I$&$-1$&  0 &  1 & $\I$&$-\I$&$-1$&  0 
\Strut\\
$(0,1)$& 20            & 1 &  0 &  1 &  0 &  0 &  0 &  1 &  0 &  1 &  0 &
                       $-1$&  0 &$-1$&  0 &  0 &  0 &$-1$&  0 &$-1$&  0 
\Strut\\
$(1,1)$& 20            & 1 &  0 & $\I$&  0 &  0 &  0 &$-\I$&  0 &$-1$&  0 &
                       $-1$&  0 &$-\I$&  0 &  0 &  0 & $\I$&  0 &  1 &  0 
\Strut\\
$(2,1)$& 4             & 1 &  0 &$-1$&  0 &  1 &  0 &$-1$&  0 &  1 &  0 &
                       $-1$&  0 &  1 &  0 &$-1$&  0 &  1 &  0 &$-1$&  0 
\Strut\\
$(3,1)$& 20            & 1 &  0 &$-\I$&  0 &  0 &  0 & $\I$&  0 &$-1$&  0 &
                       $-1$&  0 & $\I$&  0 &  0 &  0 &$-\I$&  0 &  1 &  0 
\Strut\\
\hline
\hline
       & $\ell$        & 1 &  4 &  4 &    &  1 &    &  4 &    &  2 &    &
                         2 &    &  4 &    &    &    &  4 &    &  2 &    
\Strut\\
       & $m$           & 8 &  1 &  2 &    &  2 &    &  2 &    &  4 &    &
                         4 &    &  2 &    &    &    &  2 &    &  4 &    
\Strut\\
\hline
\end{tabular*}
\end{table*}
  
Consider now a primitive Dirichlet character $\chi$. Its $L$-series is
defined as
\begin{equation}\label{lfun}
L(s,\chi)\; := \; \sum_{k=1}^{\infty}\,\frac{\chi(k)}{k^{s}}
\; =\;\sum_{k=1}^{f^{}_{\chi}}\,\chi(k)\,
\sum_{j=0}^{\infty}\,\frac{1}{(k+jf^{}_{\chi})^s}
\end{equation}
which has nice analytic properties. In particular, it converges for
all $s\!\in\!\CC$ with ${\rm Re}(s)\!>\!1$, which follows by
expressing the last sum in terms of the Hurwitz zeta function
\cite[Ch.~4]{W}.  Due to the total multiplicativity of $\chi$, its
$L$-series has a particularly simple Euler product expansion
\cite[p.~31]{W}, namely
\begin{equation}\label{leuler}
L(s,\chi)\; = \; \prod_{p}\,\frac{1}{1-\chi(p)\,p^{-s}}\, ,\quad
{\rm Re}(s)>1 \,.
\end{equation}
The following result is standard \cite[Thm.~4.3]{W}.
\begin{prop}  \label{prop-3}
The Dedekind zeta function of the cyclotomic field\/ $\QQ(\xi_{n})$ is
given by 
\[
\zeta^{}_{\QQ(\xi_{n})}(s)\;=\; \prod_{\chi\in\widehat{G}_n} L(s,\chi) 
\]
where\/ $\widehat{G}_n$ is the set of primitive Dirichlet characters of\/
$\QQ(\xi_{n})/\QQ$.
\hspace*{\fill}$\square$
\end{prop}

The Euler product expansion of $\zeta^{}_{\QQ(\xi_{n})}(s)$ is 
\begin{eqnarray}
\zeta^{}_{\QQ(\xi_{n})}(s)& =&
\prod_{\chi\in\widehat{G}_{n}}\prod_{p}\,\frac{1}{1-\chi(p)\,p^{-s}}\nonumber\\
&=&\prod_{p}\prod_{\chi\in\widehat{G}_{n}}\frac{1}{1-\chi(p)\,p^{-s}}\, ,
\end{eqnarray}
so, using Proposition~\ref{prop:dede}, our Euler factors of
Eq.~\eqref{euler} are
\begin{equation}
E_{n}(p^{-s})\;=\;
\prod_{\chi\in\widehat{G}_{n}}{\big(1-\chi(p)\,p^{-s}\big)}^{-1}\, .
\label{eulerviachar}
\end{equation}
Since $\chi$ is $f^{}_{\chi}$-periodic, and $f^{}_{\chi}|n$, all
primitive Dirichlet characters in $\widehat{G}_{n}$ are $n$-periodic,
in agreement with their original construction. Therefore, the
structure of an Euler factor can only depend on the residue class of
$p\bmod n$.

\begin{table*}
\caption[]{First terms of the Dirichlet series of Eq.~\eqref{def:diri} 
for the class number one cases as listed in \eqref{liste}.\label{diritab}}
\def\Strut{\large\strut}
\renewcommand{\arraystretch}{2.2}
\small
\begin{tabular*}{\textwidth}{@{\extracolsep{\fill}} @{} rl @{}}
\hline
 $n$ & $\zeta^{}_{{\mathcal M}_n}(s)$\Strut\\
\hline
 3 & $1 + \frac{1}{3^s} \+ \frac{1}{4^s} \+ \frac{2}{7^s} \+
      \frac{1}{9^s} \+ \frac{1}{12^s} \+ \frac{2}{13^s} \+
      \frac{1}{16^s} \+ \frac{2}{19^s} \+ \frac{2}{21^s} \+
      \frac{1}{25^s} \+ \frac{1}{27^s} \+ \frac{2}{28^s} \+
      \frac{2}{31^s} \+ \frac{1}{36^s} \+ \frac{2}{37^s} \+
      \frac{2}{39^s} \+ \frac{2}{43^s} \+ \frac{1}{48^s} \+ 
      \frac{3}{49^s} \+ \ldots$\Strut\\
 4 & $1 + \frac{1}{2^s} \+ \frac{1}{4^s} \+ \frac{2}{5^s} \+ 
      \frac{1}{8^s} \+ \frac{1}{9^s} \+ \frac{2}{10^s} \+
      \frac{2}{13^s} \+ \frac{1}{16^s} \+  \frac{2}{17^s} \+
      \frac{1}{18^s} \+ \frac{2}{20^s} \+  \frac{3}{25^s} \+ 
      \frac{2}{26^s} \+ \frac{2}{29^s} \+ \frac{1}{32^s} \+ 
      \frac{2}{34^s} \+ \frac{1}{36^s} \+ \frac{2}{37^s} \+ 
      \frac{2}{40^s} \+ \ldots$\Strut\\
 5 & $1 + \frac{1}{5^s} \+ \frac{4}{11^s} \+ \frac{1}{16^s} \+ 
      \frac{1}{25^s} \+ \frac{4}{31^s} \+ \frac{4}{41^s} \+
      \frac{4}{55^s} \+ \frac{4}{61^s} \+  \frac{4}{71^s} \+
      \frac{1}{80^s} \+ \frac{1}{81^s} \+ \frac{4}{101^s} \+ 
      \frac{10}{121^s} \+  \frac{1}{125^s} \+  \frac{4}{131^s} \+ 
      \frac{4}{151^s} \+  \frac{4}{155^s} \+ \ldots$\Strut\\
 7 & $1 + \frac{1}{7^s} \+ \frac{2}{8^s} \+ \frac{6}{29^s} \+
      \frac{6}{43^s} \+ \frac{1}{49^s} \+ \frac{2}{56^s} \+
      \frac{3}{64^s} \+ \frac{6}{71^s} \+ \frac{6}{113^s} \+
      \frac{6}{127^s} \+ \frac{3}{169^s} \+ \frac{6}{197^s} \+ 
      \frac{6}{203^s} \+ \frac{6}{211^s} \+ \frac{12}{232^s} \+ 
      \frac{6}{239^s} \+ \frac{6}{281^s} \+ \ldots $\Strut\\
 8 & $1 + \frac{1}{2^s} \+ \frac{1}{4^s} \+ \frac{1}{8^s} \+
      \frac{2}{9^s} \+ \frac{1}{16^s} \+ \frac{4}{17^s} \+
      \frac{2}{18^s} \+ \frac{2}{25^s} \+ \frac{1}{32^s} \+
      \frac{4}{34^s} \+ \frac{2}{36^s} \+ \frac{4}{41^s} \+
      \frac{2}{49^s} \+ \frac{2}{50^s} \+ \frac{1}{64^s} \+
      \frac{4}{68^s} \+ \frac{2}{72^s} \+ \frac{4}{73^s} \+ 
      \frac{3}{81^s} \+ \ldots$\Strut\\
 9 & $1 + \frac{1}{3^s} \+ \frac{1}{9^s} \+ \frac{6}{19^s} \+ 
      \frac{1}{27^s} \+ \frac{6}{37^s} \+ \frac{6}{57^s} \+ 
      \frac{1}{64^s} \+ \frac{6}{73^s} \+ \frac{1}{81^s} \+ 
      \frac{6}{109^s} \+ \frac{6}{111^s} \+ \frac{6}{127^s} \+ 
      \frac{6}{163^s} \+ \frac{6}{171^s} \+ \frac{6}{181^s} \+ 
      \frac{1}{192^s} \+ \frac{6}{199^s} \+ \ldots $\Strut\\
11 & $1 + \frac{1}{11^s} \+ \frac{10}{23^s} \+ \frac{10}{67^s} \+ 
      \frac{10}{89^s} \+ \frac{1}{121^s} \+ \frac{10}{199^s} \+ 
      \frac{2}{243^s} \+ \frac{10}{253^s} \+ \frac{10}{331^s} \+ 
      \frac{10}{353^s} \+  \frac{10}{397^s} \+ \frac{10}{419^s} \+ 
      \frac{10}{463^s} \+ \frac{55}{529^s} \+ \frac{10}{617^s} \+ 
      \frac{10}{661^s} \+ \ldots$\Strut\\
12 & $1 + \frac{1}{4^s} \+ \frac{1}{9^s} \+ \frac{4}{13^s} \+
      \frac{1}{16^s} \+ \frac{2}{25^s} \+ \frac{1}{36^s} \+
      \frac{4}{37^s} \+ \frac{2}{49^s} \+ \frac{4}{52^s} \+
      \frac{4}{61^s} \+ \frac{1}{64^s} \+ \frac{4}{73^s} \+
      \frac{1}{81^s} \+ \frac{4}{97^s} \+ \frac{2}{100^s} \+
      \frac{4}{109^s} \+ \frac{4}{117^s} \+ \frac{2}{121^s} \+
      \ldots$\Strut\\
13 & $1 + \frac{1}{13^s} \+ \frac{4}{27^s} \+ \frac{12}{53^s} \+ 
      \frac{12}{79^s} \+ \frac{12}{131^s} \+ \frac{12}{157^s}\+ 
      \frac{1}{169^s} \+ \frac{12}{313^s} \+ \frac{4}{351^s} \+ 
      \frac{12}{443^s} \+ \frac{12}{521^s} \+ \frac{12}{547^s} \+ 
      \frac{12}{599^s} \+ \frac{3}{625^s} \+ \frac{12}{677^s} \+ 
      \frac{12}{689^s} \+ \ldots$\Strut\\
15 & $1 + \frac{2}{16^s} \+ \frac{1}{25^s} \+ \frac{8}{31^s} \+ 
      \frac{8}{61^s} \+ \frac{1}{81^s} \+ \frac{4}{121^s} \+ 
      \frac{8}{151^s} \+ \frac{8}{181^s} \+ \frac{8}{211^s} \+ 
      \frac{8}{241^s} \+ \frac{3}{256^s} \+ \frac{8}{271^s} \+ 
      \frac{8}{331^s} \+ \frac{4}{361^s} \+ \frac{2}{400^s} \+ 
      \frac{8}{421^s} \+ \ldots$\Strut\\
16 & $1 + \frac{1}{2^s} \+ \frac{1}{4^s} \+ \frac{1}{8^s} \+ 
      \frac{1}{16^s} \+ \frac{8}{17^s} \+ \frac{1}{32^s} \+ 
      \frac{8}{34^s} \+ \frac{4}{49^s} \+ \frac{1}{64^s} \+
      \frac{8}{68^s} \+ \frac{2}{81^s} \+ \frac{8}{97^s} \+
      \frac{4}{98^s} \+ \frac{8}{113^s} \+ \frac{1}{128^s} \+ 
      \frac{8}{136^s} \+ \frac{2}{162^s} \+ \frac{8}{193^s} \+
      \ldots$\Strut\\
17 & $1 + \frac{1}{17^s} \+ \frac{16}{103^s} \+ \frac{16}{137^s} \+
      \frac{16}{239^s} \+ \frac{2}{256^s} \+ \frac{1}{289^s} \+
      \frac{16}{307^s} \+ \frac{16}{409^s} \+ \frac{16}{443^s} \+
      \frac{16}{613^s} \+ \frac{16}{647^s} \+ \frac{16}{919^s} \+
      \frac{16}{953^s} \+ \frac{16}{1021^s} \+ \frac{16}{1123^s} \+
      \ldots$\Strut\\
19 & $1 + \frac{1}{19^s} \+ \frac{18}{191^s} \+ \frac{18}{229^s} \+
      \frac{6}{343^s} \+ \frac{1}{361^s} \+ \frac{18}{419^s} \+
      \frac{18}{457^s} \+ \frac{18}{571^s} \+ \frac{18}{647^s} \+
      \frac{18}{761^s} \+ \frac{18}{1103^s} \+ \frac{18}{1217^s} \+
      \frac{6}{1331^s} \+ \frac{9}{1369^s} \+ \frac{18}{1483^s} \+
      \ldots$\Strut\\
20 & $1 + \frac{2}{5^s} \+ \frac{1}{16^s} \+ \frac{3}{25^s} \+ 
      \frac{8}{41^s} \+ \frac{8}{61^s} \+ \frac{2}{80^s} \+ 
      \frac{2}{81^s} \+ \frac{8}{101^s} \+ \frac{4}{121^s} \+
      \frac{4}{125^s} \+ \frac{8}{181^s} \+ \frac{16}{205^s} \+ 
      \frac{8}{241^s} \+ \frac{1}{256^s} \+ \frac{8}{281^s} \+ 
      \frac{16}{305^s} \+ \frac{4}{361^s} \+ \ldots$\Strut\\
21 & $1 + \frac{2}{7^s} \+ \frac{12}{43^s} \+ \frac{3}{49^s} \+
      \frac{2}{64^s} \+ \frac{12}{127^s} \+ \frac{6}{169^s}\+
      \frac{12}{211^s} \+ \frac{24}{301^s} \+ \frac{12}{337^s} \+ 
      \frac{4}{343^s} \+ \frac{12}{379^s} \+ \frac{12}{421^s} \+ 
      \frac{4}{448^s} \+ \frac{12}{463^s} \+ \frac{12}{547^s} \+ 
      \frac{12}{631^s} \+ \ldots$\Strut\\
24 & $1 + \frac{1}{4^s} \+ \frac{2}{9^s} \+ \frac{1}{16^s} \+ 
      \frac{4}{25^s} \+ \frac{2}{36^s} \+ \frac{4}{49^s} \+ 
      \frac{1}{64^s} \+ \frac{8}{73^s} \+ \frac{3}{81^s} \+
      \frac{8}{97^s} \+ \frac{4}{100^s} \+ \frac{4}{121^s} \+ 
      \frac{2}{144^s} \+\frac{4}{169^s} \+ \frac{8}{193^s} \+ 
      \frac{4}{196^s} \+ \frac{8}{225^s} \+ \ldots$\Strut\\
25 & $1 + \frac{1}{5^s} \+ \frac{1}{25^s} \+ \frac{20}{101^s} \+
      \frac{1}{125^s} \+ \frac{20}{151^s} \+ \frac{20}{251^s} \+
      \frac{20}{401^s} \+ \frac{20}{505^s} \+ \frac{20}{601^s} \+
      \frac{1}{625^s} \+ \frac{20}{701^s} \+ \frac{20}{751^s} \+
      \frac{20}{755^s} \+ \frac{20}{1051^s} \+ \frac{20}{1151^s} \+
      \frac{20}{1201^s} \+ \ldots$\Strut\\
27 & $1 + \frac{1}{3^s} \+ \frac{1}{9^s} \+ \frac{1}{27^s} \+
      \frac{1}{81^s} \+ \frac{18}{109^s} \+ \frac{18}{163^s} \+
      \frac{1}{243^s} \+ \frac{18}{271^s} \+ \frac{18}{327^s} \+
      \frac{18}{379^s} \+ \frac{18}{433^s} \+ \frac{18}{487^s} \+
      \frac{18}{489^s} \+ \frac{18}{541^s} \+ \frac{1}{729^s} \+
      \frac{18}{757^s} \+ \ldots$\Strut\\
28 & $1 + \frac{2}{8^s} \+ \frac{12}{29^s} \+ \frac{1}{49^s} \+
      \frac{3}{64^s} \+ \frac{12}{113^s} \+ \frac{6}{169^s}\+ 
      \frac{12}{197^s} \+ \frac{24}{232^s} \+ \frac{12}{281^s} \+ 
      \frac{12}{337^s} \+ \frac{2}{392^s} \+ \frac{12}{421^s} \+ 
      \frac{12}{449^s} \+ \frac{4}{512^s} \+ \frac{12}{617^s} \+ 
      \frac{12}{673^s} \+ \ldots$\Strut\\
32 & $1 + \frac{1}{2^s} \+ \frac{1}{4^s} \+ \frac{1}{8^s} \+
      \frac{1}{16^s} \+ \frac{1}{32^s} \+ \frac{1}{64^s} \+
      \frac{16}{97^s} \+ \frac{1}{128^s} \+ \frac{16}{193^s} \+
      \frac{16}{194^s} \+ \frac{1}{256^s} \+ \frac{16}{257^s} \+
      \frac{8}{289^s} \+ \frac{16}{353^s} \+ \frac{16}{386^s} \+
      \frac{16}{388^s} \+ \frac{16}{449^s} \+ \ldots$\Strut\\
33 & $1 + \frac{20}{67^s} \+ \frac{1}{121^s} \+ \frac{20}{199^s} \+
      \frac{2}{243^s} \+ \frac{20}{331^s} \+ \frac{20}{397^s} \+
      \frac{20}{463^s} \+ \frac{10}{529^s} \+ \frac{20}{661^s} \+
      \frac{20}{727^s} \+ \frac{20}{859^s} \+ \frac{20}{991^s} \+
      \frac{2}{1024^s} \+ \frac{20}{1123^s} \+ \frac{20}{1321^s} \+
      \ldots$\Strut\\
35 & $1 + \frac{24}{71^s} \+ \frac{24}{211^s} \+ \frac{24}{281^s} \+
      \frac{24}{421^s} \+ \frac{24}{491^s} \+ \frac{24}{631^s} \+
      \frac{24}{701^s} \+ \frac{12}{841^s} \+ \frac{24}{911^s} \+
      \frac{24}{1051^s} \+ \frac{8}{1331^s} \+ \frac{24}{1471^s} \+
      \frac{12}{1681^s} \+ \frac{24}{2311^s} \+ \frac{24}{2381^s} \+
      \ldots$\Strut\\
36 & $1 + \frac{1}{9^s} \+ \frac{12}{37^s} \+ \frac{1}{64^s} \+ 
      \frac{12}{73^s} \+ \frac{1}{81^s} \+ \frac{12}{109^s}\+
      \frac{12}{181^s} \+ \frac{6}{289^s} \+ \frac{12}{333^s} \+ 
      \frac{6}{361^s} \+ \frac{12}{397^s} \+ \frac{12}{433^s} \+ 
      \frac{12}{541^s} \+ \frac{1}{576^s} \+ \frac{12}{577^s} \+ 
      \frac{12}{613^s} \+ \ldots$\Strut\\
40 & $1 + \frac{1}{16^s} \+ \frac{2}{25^s} \+ \frac{16}{41^s} \+
      \frac{4}{81^s} \+ \frac{8}{121^s} \+ \frac{16}{241^s} \+
      \frac{1}{256^s} \+ \frac{16}{281^s} \+ \frac{8}{361^s} \+
      \frac{2}{400^s} \+ \frac{16}{401^s} \+ \frac{16}{521^s} \+
      \frac{16}{601^s} \+ \frac{3}{625^s} \+ \frac{16}{641^s} \+
      \frac{16}{656^s} \+ \ldots$\Strut\\
44 & $1 + \frac{20}{89^s} \+ \frac{1}{121^s} \+ \frac{20}{353^s} \+
      \frac{20}{397^s} \+ \frac{10}{529^s} \+ \frac{20}{617^s} \+
      \frac{20}{661^s} \+ \frac{20}{881^s} \+ \frac{20}{1013^s} \+
      \frac{1}{1024^s} \+ \frac{20}{1277^s} \+ \frac{20}{1321^s} \+
      \frac{20}{1409^s} \+ \frac{20}{1453^s} \+ \frac{10}{1849^s} \+
      \ldots$\Strut\\
45 & $1 + \frac{1}{81^s} \+ \frac{24}{181^s} \+ \frac{24}{271^s} \+
      \frac{12}{361^s} \+ \frac{24}{541^s} \+ \frac{24}{631^s} \+
      \frac{24}{811^s} \+ \frac{24}{991^s} \+ \frac{24}{1171^s} \+
      \frac{24}{1531^s} \+ \frac{24}{1621^s} \+ \frac{24}{1801^s} \+
      \frac{24}{2161^s} \+ \frac{24}{2251^s} \+ \frac{24}{2341^s} \+
      \ldots$\Strut\\
48 & $1 + \frac{1}{4^s} \+ \frac{1}{16^s} \+ \frac{8}{49^s} \+
      \frac{1}{64^s} \+ \frac{2}{81^s} \+ \frac{16}{97^s} \+
      \frac{16}{193^s} \+ \frac{8}{196^s} \+ \frac{16}{241^s} \+
      \frac{1}{256^s} \+ \frac{8}{289^s} \+ \frac{2}{324^s} \+
      \frac{16}{337^s} \+ \frac{16}{388^s} \+ \frac{16}{433^s} \+
      \frac{8}{529^s} \+ 
      \ldots$\Strut\\
60 & $1 + \frac{2}{16^s} \+ \frac{2}{25^s} \+ \frac{16}{61^s} \+
      \frac{2}{81^s} \+ \frac{8}{121^s} \+ \frac{16}{181^s} \+
      \frac{16}{241^s} \+ \frac{3}{256^s} \+ \frac{8}{361^s} \+
      \frac{4}{400^s} \+ \frac{16}{421^s} \+ \frac{16}{541^s} \+
      \frac{16}{601^s} \+ \frac{3}{625^s} \+ \frac{16}{661^s} \+
      \frac{8}{841^s} \+ \ldots$\Strut\\
84 & $1 + \frac{2}{49^s} \+ \frac{2}{64^s} \+ \frac{12}{169^s} \+
      \frac{24}{337^s} \+ \frac{24}{421^s} \+ \frac{24}{673^s} \+
      \frac{2}{729^s} \+ \frac{24}{757^s} \+ \frac{12}{841^s} \+
      \frac{24}{1009^s} \+ \frac{24}{1093^s} \+ \frac{24}{1429^s} \+
      \frac{24}{1597^s} \+ \frac{12}{1681^s} \+ \frac{12}{1849^s} \+
      \ldots$\Strut\\[2mm] 
\hline
\end{tabular*}
\renewcommand{\arraystretch}{1}
\end{table*}

Calculating the Euler factors $E_{n}(p^{-s})$, which means to expand
the product on the right-hand side of Eq.~\eqref{eulerviachar}, one
finds that they are always of the form ${(1-p^{-\ell s})}^{-m}$ with
integers $\ell$ and $m$, compare Eq.~\eqref{factor}. These integers,
which are now simple arithmetic expressions in the values of the
characters at $p$, are precisely the quantities introduced after
Eq.~\eqref{factor}. This works for all primes, and the distinction
between $p|n$ and $p\!\nmid\! n$ is implicit, see Table~\ref{tab:char}
for an example.

Whenever a prime $p$ is not ramified, which happens if and only if
$\chi(p)\ne 0$ for all $\chi\in\widehat{G}_{n}$, these integers, as
mentioned above, satisfy $\ell\, m=\phi(n)$. For the remaining primes,
which are precisely those dividing $n$, the product $\ell\, m$ is a
true divisor of $\phi(n)$, and counts the number of characters not
vanishing at $p$.  In the example of Table~\ref{tab:char}, these are
the primes $p=2$ and $p=5$, while all other primes are
unramified. Further details, together with the meaning of $\ell$ and
$m$ for the splitting of primes in the algebraic field extension, can
be found in \cite{L,W}. An explicit calculation can easily be done
with an algebraic program package; a corresponding
Mathematica$^{\circledR}$ program for the calculation of the
generating functions can be downloaded from \cite{u}.

We summarise the result as follows.

\begin{thm}
The Euler factors of the Dedekind zeta function of\/ $\QQ(\xi_{n})$
are of the form given in Eq.~$\eqref{factor}$ with characteristic
indices\/ $\ell$ and\/ $m$ that depend on\/ $p$ and\/ $n$. Using the
explicit representation via the\/ $L$-series, the results for the
class number one cases are those given in Tables~$\ref{tab:basic}$
and\/ $\ref{tab:rami}$.\hspace*{\fill}$\square$
\end{thm}

Since these results have concrete applications in crystallography and
materials science, we spell out the first few terms of the generating
functions in Table~\ref{diritab}. Explicit realisations of the
corresponding colourings can be constructed by means of prime (ideal)
factorisation in $\ZZ[\xi_n]$ and interpretation of the ideals
as similarity submodules.

\section{Asymptotic properties}

Generating functions are an efficient way to encapsulate an entire
series of numbers into one object. If the generating function is not
only a formal object, but also has well-defined analytic properties,
one can extract asymptotic properties of the sequence from the
generating function. In the case of Dirichlet series, this requires
some relatively advanced techniques from complex
analysis. Fortunately, in the cases that emerged here, we are still in
a rather simple situation, fully covered by the following special case
of Delange's theorem, see \cite[Ch.\ II.7, Thm.\ 15]{T}.
\begin{prop}\label{delange}
Let\/ $F(s)=\sum_{k=1}^{\infty}a(k)\,k^{-s}$ be a Dirichlet series
with non-negative coefficients which converges for\/ ${\rm
Re}(s)>1$. Suppose that\/ $F(s)$ is holomorphic at all points of the
line\/ $\{ {\rm Re}(s)=1\}$ except at\/ $s=1$. Approaching\/ $s=1$
from the half-plane\/ $\{ {\rm Re}(s)>1\}$, let\/ $F(s)$ have a
singularity of the form\/ $F(s) = g(s) + h(s)/(s-1)$, where both\/
$g(s)$ and\/ $h(s)$ are holomorphic at\/ $s=1$. Then, as\/
$x\rightarrow\infty$, 
\[
     A(x) \; := \; \sum_{k\leq x} a(k)
          \; \sim \;  h(1)\, x \, .
\]
In other words, $h(1)$ is the asymptotic average value of\/ $a(k)$
as\/ $k\to\infty$.
\hspace*{\fill}$\square$
\end{prop}

Our Dirichlet series $F_{n}(s)$ are Dedekind zeta functions of
cyclotomic fields, according to Proposition~\ref{prop:dede}.  Then,
Proposition~\ref{delange} leads to the following explicit result.

\begin{thm}\label{thm:asym}
Let\/ $n$ be a number from the list~$\eqref{liste}$. The average number
of Bravais colourings of\/ ${\mathcal{M}}_{n}$, with a given number\/
$k$ of colours, is asymptotically, as\/ $k\to\infty$, given by the
residue of the Dedekind zeta function of\/ $\QQ(\xi_{n})$ at\/
$s=1$. Consequently, for\/  $x\to\infty$, one has
\[
A_{n}(x)\;:=\;\sum_{k\leq x}
a^{}_{n}(k)\;\sim\;\alpha^{}_{n}\, x\, ,
\]
with
\[
\alpha^{}_{n}\;=\; \prod_{1\not\equiv\chi\in\widehat{G}_{n}}L(1,\chi)\, .
\]
Some values of\/ $\alpha^{}_{n}$ are given in Table~$\ref{tab:asym}$.
\end{thm}
\pf
The Dedekind zeta functions of the cyclotomic fields (resp.\ their
analytic continuations) always possess an isolated pole of first order
at $s=1$, and no other singularity in the entire complex plane, see
\cite[Ch.~4]{W}. Consequently, by Proposition~\ref{delange},
$A_{n}(x)$ has linear growth as $x\to\infty$. The growth rate is then
the residue of the zeta function at $s=1$. The $L$-series $L(s,\chi)$
have analytic continuations to entire functions for all Dirichlet
characters $\chi$ except for $\chi^{}_{0}\equiv 1$. In the latter
case, one has $L(s,\chi^{}_{0})=\zeta(s)$. This is Riemann's zeta
function which has a simple pole at $s=1$ with residue $1$. For all
$\chi\ne\chi_{0}$, one has $L(1,\chi)\ne 0$, see \cite[Cor.~4.4]{W},
and the expression for $\alpha^{}_{n}$ is then immediate from
Proposition~\ref{prop-3}.
\qed

\begin{table}
\caption[]{Residues of Dedekind zeta functions.\label{tab:asym}}\smallskip
\def\Strut{\large\strut}
\small
\renewcommand{\arraystretch}{1.75}
\begin{tabular*}{\linewidth}{@{\extracolsep{\fill}} @{} 
cr@{\qquad}c@{\qquad}c@{\qquad}l @{}}
\hline
$\phi(n)$ & $n$ & $R_{n}$ & residue & numerical \Strut\\
\hline
2 &  3 & $1$ & $\frac{\pi\sqrt{3}}{3^2}$ & 
         $0.604\,600$\rule[0ex]{0ex}{3.5ex}\Strut\\
  &  4 & $1$ & $\frac{\pi}{2^2}$  &  
$0.785\,398$\rule[-2ex]{0ex}{2ex}\Strut\\
\hline
4 &  5 & $2\log(\frac{1+\sqrt{5}}{2})$ &
         $\frac{2\pi^{2}\sqrt{5}}{5^3}R_{4}$ & 
         $0.339\,837$\rule[0ex]{0ex}{3.5ex}\Strut\\
  &  8 & $2\log(1+\sqrt{2\,})$ & 
         $\frac{\pi^{2}}{2^5} R_{8}$ & $0.543\,676$\Strut\\
  & 12 & $\log(2+\sqrt{3\,})$ & 
         $\frac{\pi^{2}}{2^{2}3^{2}}R_{12}$ & 
         $0.361\,051$\rule[-2ex]{0ex}{2ex}\Strut\\
\hline
6 &  7 & $\protect\hphantom{11}2.101\,819$ & 
         $\frac{2^2\pi^{3}\sqrt{7}}{7^4} R_{7}$ & 
         $0.287\,251$\rule[0ex]{0ex}{3.5ex}\Strut\\
  &  9 & $\protect\hphantom{11}3.397\,150$& 
         $\frac{2^2\pi^{3}\sqrt{3}}{3^7} R_{9}$ & 
         $0.333\,685$\rule[-2ex]{0ex}{2ex}\Strut\\
\hline
8 & 15 & $\protect\hphantom{11}4.661\,821$  & 
         $\frac{2^3\pi^4}{3^{3}5^{4}}R_{15}$ & 
         $0.215\,279$\rule[0ex]{0ex}{3.5ex}\Strut\\
  & 16 & $\protect\hphantom{1}19.534\,360$ & 
         $\frac{\pi^{4}}{2^{12}}R_{16}$ & $0.464\,557$\Strut\\
  & 20 & $\protect\hphantom{11}7.411\,242$  & 
         $\frac{\pi^4}{2^{2}5^{4}}R_{20}$ & $0.288\,769$\Strut\\
  & 24 & $\protect\hphantom{1}10.643\,594$ & 
         $\frac{\pi^{4}}{2^{7}3^{3}}R_{24}$ & 
         $0.299\,995$\rule[-2ex]{0ex}{2ex}\Strut\\
\hline
10& 11 & $\protect\hphantom{1}26.171\,106$ &
         $\frac{2^{4}\pi^{5}\sqrt{11}}{11^{6}}R_{11}$ & 
         $0.239\,901$\rule[-2ex]{0ex}{5.5ex}\Strut\\
\hline
12& 13 & $120.784\,031$  & 
         $\frac{2^{5}\pi^{6}\sqrt{13}}{13^{7}}R_{13}$ & 
         $0.213\,514$\rule[0ex]{0ex}{3.5ex}\Strut\\
  & 21 & $\protect\hphantom{1}70.399\,398$ & 
         $\frac{2^{5}\pi^{6}}{3^{4}7^{6}}R_{21}$ & 
         $0.227\,271$\Strut\\
  & 28 & $123.252\,732$  & 
         $\frac{\pi^{6}}{2^{2}7^{6}}R_{28}$ & 
         $0.251\,795$\Strut\\
  & 36 & $162.837\,701$ & 
         $\frac{\pi^{6}}{2^{2}3^{11}}R_{36}$ & 
         $0.220\,933$\rule[-2ex]{0ex}{2ex}\Strut\\
\hline
\end{tabular*}
\renewcommand{\arraystretch}{1}
\end{table}

To apply this Theorem, we need to know how to calculate the values
$L(1,\chi)$ for the non-trivial characters. This is described
explicitly on p.~36 and in Thm.~4.9 of \cite{W}. The result is
\begin{equation}
\alpha^{}_{n}\;=\;\frac{R_{n}}{N(n)}\,
\left(\frac{2\pi\prod_{p|n}p^{1/(p-1)}}{n}\right)^{\phi(n)/2}\, ,
\end{equation}
where $p$ in the product runs over all prime divisors of $n$, hence
precisely over the ramified primes, $N(n)$ is the function of
Eq.~\eqref{symm}, and $R_{n}$ is the {\em regulator}\/ of the
cyclotomic field $\QQ(\xi_{n})$. Its calculation, which is a bit
technical, is explained in \cite{W}. It is based on the knowledge of
the fundamental units of the maximal order of
$\QQ(\xi_{n}+\bar{\xi}_{n})$, the maximal real subfield of
$\QQ(\xi_{n})$, followed by the calculation of a determinant
\cite[p.~41 and Prop.~4.16]{W}.  We show the first few cases in
Table~\ref{tab:asym}.  For $6\le \phi(n)\le 12$, the fundamental units
needed for the (numerical) calculation of the regulator were
determined by means of the program package KANT \cite{K1,K2}.

\begin{acknowledgment}
It is a pleasure to thank Robert V.\ Moody and Peter A.\ B.\ Pleasants
for cooperation and helpful discussions. We express our gratitude to
the Erwin Schr\"{o}dinger International Institute for Mathematical
Physics in Vienna for support during a stay in 2002/2003, where most
of this work was done.
\end{acknowledgment}


\begin{thebibliography}{99}\itemsep 0pt

\bibitem{B}
M.~Baake,
Combinatorial aspects of colour symmetries,
{\em J.\ Phys.\ A: Math.\ Gen.}\/ {\bf 30} (1997) 2687--2698;
{\tt mp\_arc/02-323}.

\bibitem{BG}
M.~Baake and U.~Grimm,
Combinatorial problems of (quasi)\-crystallography,
in: {\em Quasicrystals: Structure and Physical
Properties}, ed.\ H.-R.\ Trebin
(Wiley-VCH, Weinheim, 2003) pp.~160--171;
{\tt math-ph/0212015}.

\bibitem{BGS1}
M.~Baake, U.~Grimm and M.~Scheffer,
Colourings of planar quasicrystals,
{\em J.\ Alloys Comp.}\/ {\bf 342} (2002) 195--197;
{\tt cond-mat/0110654}.

\bibitem{BGS2}
M.~Baake, U.~Grimm and M.~Scheffer,
Tilings colourful --- only even more so,
{\em J.\ Non-Cryst.\ Solids}\/ {\bf 334 \& 335}
(2004) 83--85; {\tt mp\_arc/03-201}.

\bibitem{BHLS}
M.~Baake, J.~Hermisson, R.~L\"{u}ck and M.~Scheffer,
Colourings of quasicrystals, in: 
{\em Quasicrystals}, eds.\ S.~Takeuchi and T.~Fujiwara
(World Scientific, Singapore, 1998) pp.~120--123.

\bibitem{BM}
M.~Baake and R.~V.~Moody,
Similarity submodules and semigroups, in:
{\em Quasicrystals and Discrete Geometry}, ed.\ J.~Patera,
Fields Institute Monographs, vol.~10 (AMS, Providence, RI 1998)
pp.~1--13.

\bibitem{BM2}
M.~Baake and R.~V.~Moody,
Similarity submodules and root systems in four dimensions,
{\em Can.\ J.\ Math.}\/ {\bf 51} (1999) 1258--1276; 
{\tt math.MG/9904028}.

\bibitem{BS} Z.~I.~Borevich and I.~R.~Shafarevich,
{\em Number Theory}\/ (Academic Press, New York, 1966).

\bibitem{K1}
M.~Daberkow, C.~Fieker, J.~Kl\"{u}ners, M.~Pohst, K.~Roegner and 
K.~Wildanger,
KANT V4,
{\em J.\ Symbolic Comp.}\/ {\bf 24} (1997) 267--283. 

\bibitem{u}
U.~Grimm, {\sf http://mcs.open.ac.uk/ugg2/colsymm/}.

\bibitem{K2}
KANT/KASH, {\sf http://www.math.tu-berlin.de/\~{}kant/kash.html}.

\bibitem{L}
S.~Lang,
{\em Algebraic Number Theory}, 2nd ed., Springer, New York (1994);
corr.\ 3rd printing (2000).

\bibitem{L1}
R.~Lifshitz,
Theory of color symmetry for periodic and quasiperiodic crystals,
{\em Rev.\ Mod.\ Phys.}\/ {\bf 69} (1997) 1181--1218.

\bibitem{L2}
R.~Lifshitz,
Lattice color groups of quasicrystals, in: 
{\em Quasicrystals}, eds.\ S.~Takeuchi and T.~Fujiwara
(World Scientific, Singapore, 1998) pp.~103--107.

\bibitem{MRW}
N.~D.~Mermin, D.~S.~Rokhsar and D.~C.~Wright,
Beware of $46$-fold symmetry: The classification of two-dimensional
quasicrystallographic lattices,
{\em Phys.\ Rev.\ Lett.}\/ {\bf 58} (1987) 2099--2101.

\bibitem{MP}
R.~V.~Moody and J.~Patera,
Coloring of quasi-crystals,
{\em Can.\ J.\ Phys.}\/ {\bf 72} (1994) 442--452.

\bibitem{P}
P.~A.~B.~Pleasants, 
Designer quasicrystals:\thinspace Cut-and-project sets
with pre-assigned properties,
in:  {\em Directions in Mathematical Quasicrystals},
eds.\ M.\ Baake and R.\ V.\ Moody,
CRM Monograph Series, vol.~13 
(AMS, Providence, RI, 2000) pp.~95--141.

\bibitem{PBR}
P.~A.~B.~Pleasants, M.~Baake and J.~Roth,
Planar coincidences for $N$-fold symmetry,
{\em J.\ Math.\ Phys.}\/ {\bf 37} (1996) 1029--1058;
rev.\ version, {\tt math.MG/0511147}.

\bibitem{R}
J.~S.~Rutherford,
The enumeration and symmetry-significant properties of derivative lattices,
{\em Acta Cryst.\ A}\/ {\bf 48} (1992) 500--508;
2.~Classification by color lattice group,
{\em Acta Cryst.\ A}\/ {\bf 49} (1993) 293--300;
3.~Periodic colorings of a lattice, 
{\em Acta Cryst.\ A}\/ {\bf 51} (1995) 672--679.

\bibitem{SL}
M.~Scheffer and R.~L\"{u}ck,
Coloured quasiperiodic patters with tenfold symmetry and eleven colours,
{\em J.\ Non-Cryst.\ Solids}\/ {\bf 250--252} (1999) 815--819.

\bibitem{S1}
R.~L.~E.~Schwarzenberger, 
{\em $N$-dimensional Crystallography}\/
(Pitman, San Francisco, 1980).

\bibitem{S2}
R.~L.~E.~Schwarzenberger, 
Colour symmetry, 
{\em Bull.\ London Math.\ Soc.}\/ {\bf 16} (1984) 209--240.

\bibitem{Marj1}
M.~Senechal,
Point groups and color symmetry,
{\em Z.\ Kristallogr.}\/ {\bf 172} (1975) 1--23.

\bibitem{Marj2}
M.~Senechal,
Color groups,
{\em Discrete Appl.\ Math.}\/ {\bf 1} (1979) 51--73.

\bibitem{T}
G.~Tenenbaum, 
{\em Introduction to Analytic and Probabilistic Number Theory}\/ 
(Cambridge University Press, Cambridge, 1995).

\bibitem{W}
L.~C.~Washington, 
{\em Introduction to Cyclotomic Fields}, 2nd ed.\
(Springer, New York, 1997).

\bibitem{Wilf}
H.~Wilf,
{\em Generatingfunctionology},  2nd ed.\ 
(Academic Press, Boston MA, 1994).

\end{thebibliography}
\end{document}